\documentclass[10pt,twocolumn]{article}

\usepackage{shock}
\usepackage{graphicx}
\usepackage{latexsym}
\usepackage{layout}
\usepackage{amssymb}
\usepackage{amsmath}

\input amssym.def
\input amssym.tex

\newcounter{th}
\newenvironment{theorem}{\par \refstepcounter{th} {\sc
Theorem \arabic{th}.}\it}

\newcounter{df}
\newenvironment{defin}{\par \refstepcounter{df} {\sc
Definition.}}

\newcounter{nt}
\newenvironment{remark}{\par \refstepcounter{nt} {\sc
Remark.}}

\newcounter{lm}
\newenvironment{lemma}{\par \refstepcounter{lm} {\sc
Lemma \arabic{lm}.}\it}

\newenvironment{proof}{\par \refstepcounter{nt} {\sc
 Proof.}}{$\Box$}

\newcommand{\R}{\mathbb{R}}
\newcommand{\C}{\mathbb{C}}

\title{Singularities of convex hulls of smooth hypersurfaces}

\author{
{\sc Ilya A. Bogaevsky}\\ \\
{\small Independent University of Moscow}\\
{\small Bolsho{\u\i} Vlas$'$evski{\u\i} per.\,11, Moscow 121002, Russia}\\
{\small \it E-mail: bogaevsk@mccme.ru}}

\date{}

\begin{document}

\maketitle

\begin{abstract} \noindent
We describe singularities of the convex hull of a generic compact
smooth hypersurface in four-dimensional affine space up to
diffeomorphisms. It turns out there are only two new singularities
(in comparison with the previous dimension case) which appear at
separate points of the boundary of the convex hull and are not
removed by a small perturbation of the original hypersurface. The
first singularity does not contain functional, but has at least
nine continuous number invariants. A normal form which does not
contain invariants at all is found for the second singularity.

\end{abstract}\

\section*{Introduction}

\begin{figure}[h]
\begin{center}
\includegraphics{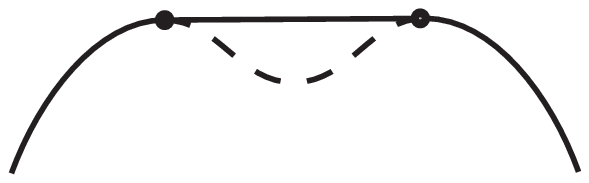}
\caption{Singularities of the convex hull of a plane curve}
\label{f1}
\end{center}
\end{figure}

The {\it convex hull\/} of a compact subset of an affine space is the
intersection of all closed half-spaces which contain the subset. The
boundary of the convex hull of a compact smooth hypersurface can have
singularities. For example, the singularities of the boundary of the
convex hull of a generic closed plane curve are discontinuities of the
second derivative (Fig.~\ref{f1}).

In the present paper we describe, up to diffeomorphism,
singularities of the convex hulls of generic compact
$C^\infty$-smooth hypersurfaces without boundary embedded into
four-dimensional affine space. A {\it singularity\/} of a convex
hull is its germ at a singular point of the boundary. As usual,
{\it generic\/} hypersurfaces are embeddings which form an open
everywhere dense set in the $C^\infty$-space of all embeddings
considered. In other words, we are interested only in
singularities which are not removed by any $C^\infty$-small
perturbation of the original hypersurface.

It turns out there are only two new singularities (in comparison
with dimension three) which appear at isolated points of the
boundary of the convex hull and are not removed by small
perturbation of the original hypersurface. The first singularity
does not contain functional moduli, but has at least nine
numerical ones. A normal form which does not contain moduli at all
is found for the second singularity.

Moreover, we show that the boundary of the convex hull is the front of a
Legendre variety and we find normal forms of its germs with respect to
contact diffeomorphisms. All singularities of the Legendre variety
prove to be stable and simple in contrast to the singularities of the
convex hull itself.

A tangent hyperplane is called a {\it support hyperplane\/} if the
hypersurface lies entirely in one of the two closed half-spaces
defined by the hyperplane. A support hyperplane is called {\it
non-singular\/} if it has only one common point with the
hypersurface and this point is a point of non-degenerate tangency.
All the remaining support planes are called {\it singular\/}. For
example, the singular support hyperplanes to a generic plane curve
are straight lines of double non-degenerate tangency.
Singularities of the convex hull of a compact hypersurface can
appear only in its singular support hyperplanes.

\medskip

\begin{figure}[h]
\begin{center}
\includegraphics{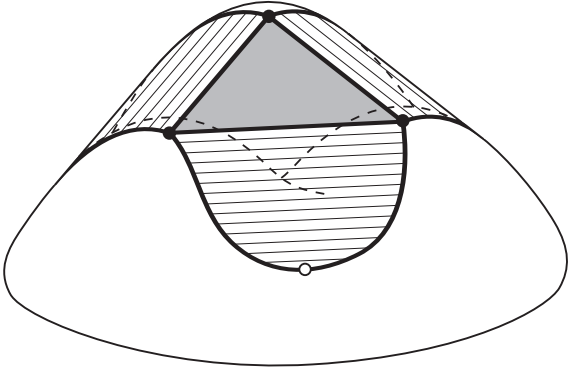}
\caption{Singularities of the convex hull of a surface} \label{f2}
\end{center}
\end{figure}

{\bf Three-dimensional space.} The convex hull of a generic
compact surface in three-dimensional space can have only two kinds
of singularities which we call {\it simplest\/} and {\it angle
singularities\/}. The normal form of the angle singularity
contains a numerical {\it modulus\/} (continuous invariant) with
respect to diffeomorphisms. These singularities shown in
Fig.~\ref{f2} are found in \cite{Z}. The results of this paper are
the following.

A typical singular support hyperplane to a generic surface in
three-dimensional space has two (and only two) points of
non-degenerate tangency with the surface ({\it $2A_1$-plane\/}).
Isolated singular support hyperplanes can have either three points
of non-degenerate tangency with the surface which form a triangle
or one point of degenerate tangency $A_3$. Such hyperplanes are
called {\it $3A_1$-} and {\it $A_3$-planes\/} respectively.

The segment between the points of tangency of a surface and its
support $2A_1$-plane is called a {\it support segment\/}. It lies
entirely in the boundary of the convex hull of the original
surface. The boundary is smooth at the interior points of the
support segment. In neighborhoods of the endpoints of the support
segment the boundary is diffeomorphic to the product of a line and
a curve with the singularity shown in Fig.~\ref{f1}. Such
singularities of a convex hull are called {\it simplest\/} and
denoted by ${\cal R}_1$.

The triangle with vertices at the points of tangency of a surface and its
support $3A_1$-plane is called a {\it support triangle\/}. In its neighborhood the
boundary of the convex hull of the surface consists of the support
triangle itself, three smooth surfaces webbed from support segments and
adjoining the sides of the support triangle, and three parts of the
original surface adjoining the vertices of the support triangle. The
convex hull has the simplest singularities at the interior points of the
sides of the support triangle. To describe singularities at its vertices
let us define the {\it under-graph\/} as the set of points lying under the graph
of a given function. It turns out that in a neighborhood of each vertex of
the support triangle the convex hull is diffeomorphic to the under-graph
of the square of the distance to an angle of a value $\beta$ where
$0<\beta<\pi$ is a unique modulus in the normal form. Such singularities
of a convex hull are called {\it angle singularities\/}. Let us note that we can consider
the epigraph of the square of the distance instead of its under-graph ---
they become diffeomorphic after the substitution $\beta \mapsto
\pi-\beta$.

Finally, the convex hull of a generic surface has the simplest singularity
at any {\it support\/} $A_3$-point that is a point of tangency of the
surface and its support $A_3$-plane. In a neighborhood of such a point
the boundary of the convex hull consists of a part of the original surface
and a surface webbed from support segments which degenerate into the
support $A_3$-point.

\medskip

{\bf Four-dimensional space.} Some singularities of the convex
hull of a generic hypersurface in four-dimensional space are
investigated in \cite{S1} and \cite{S3}. In \cite{S1} normal forms
of singularities of convex hulls are found in the case when the
support hyperplane is tangent to the original hypersurface in one
of the ways described above for three-dimensional space. Thus the
simplest and angle singularities appear in four-dimensional space
as well.

The boundary of the convex hull is smooth at the interior points of the
support segments of the $2A_1$-planes and has the simplest singularities
at their endpoints. The support segments themselves lie entirely in the
boundary of the convex hull of the original hypersurface.

At the interior points of the support triangles of the
$3A_1$-planes the boundary of the convex hull is smooth and has
the simplest singularities at the interior points of their sides.
The support triangles themselves lie entirely in the boundary of
the convex hull of the original hypersurface. In a neighborhood of
each of their vertices the convex hull is diffeomorphic to the
curvilinear cylinder over the above angle singularity with
$\beta(z)=\beta_0+z$, $\beta(z)=\beta_0+z^2$, or
$\beta(z)=\beta_0-z^2$ where $0<\beta_0<\pi$ is a unique modulus
in each of the three normal forms and $z$ is a coordinate along
the element of the cylinder. Such singularities of a convex hull
are called {\it angle singularities\/} and denoted by ${\cal
R}_2^0$, ${\cal R}_2^+$, and ${\cal R}_2^-$ respectively.

Finally, the convex hull of a generic hypersurface has the simplest
singularity at any support $A_3$-point again.

However, among the support hyperplanes of a three-dimensional hypersurface
there can be new {\it $4A_1$-} and {\it $A_1A_3$-planes\/} which are
not removed by any small perturbation of the hypersurface. Each
$4A_1$-plane has four points of non-degenerate tangency with the
hypersurface which form a tetrahedron. Each $A_1A_3$-plane has one point
of non-degenerate tangency and one point of tangency $A_3$.

\begin{figure*}[t]
\begin{center}
\includegraphics{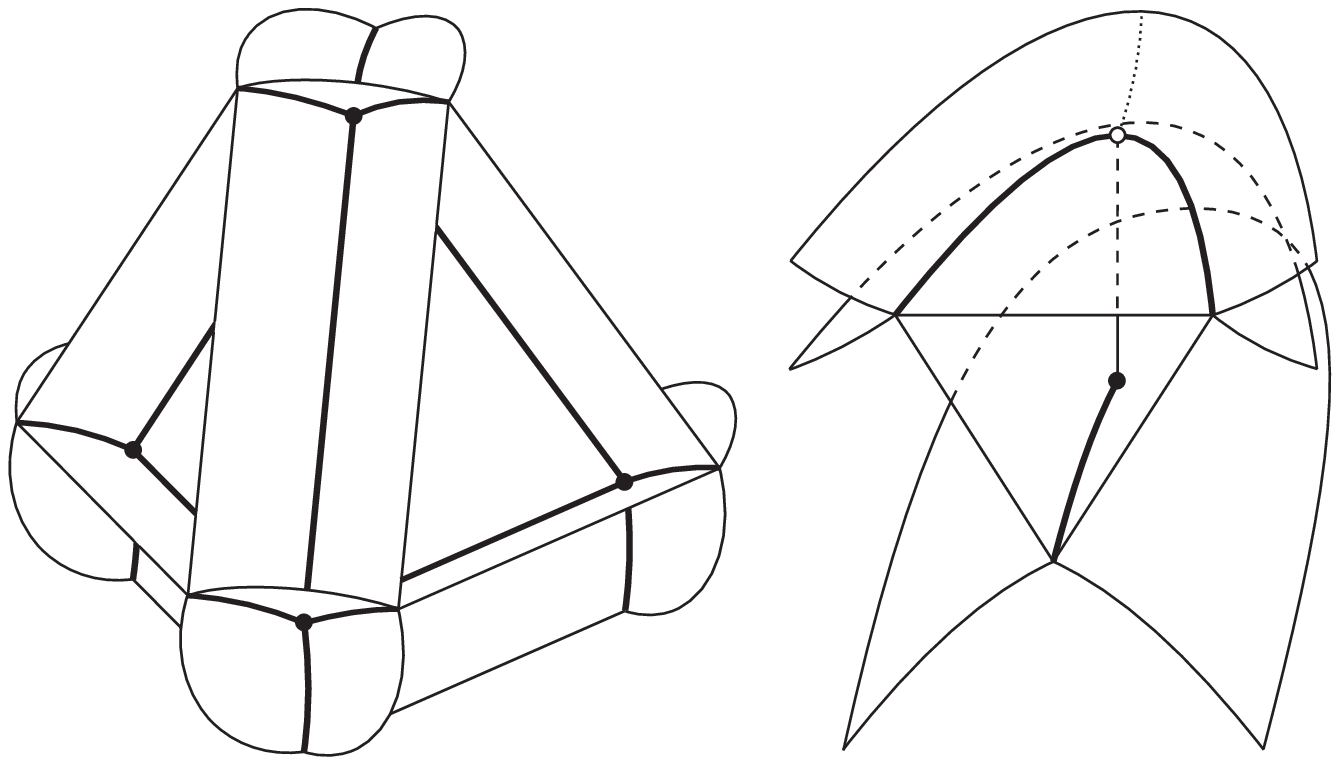}
\caption{Singularities of the convex hull of a three-dimensional
hypersurface} \label{f3}
\end{center}
\end{figure*}

The tetrahedron with vertices at the points of tangency of a
hypersurface and its support $4A_1$-plane is called a {\it support
tetrahedron\/}. In its neighborhood the boundary of the convex
hull of the hypersurface consists of the support tetrahedron
itself, four smooth strata webbed from support triangles and
adjoining the faces of the support tetrahedron, six smooth strata
webbed from support segments and adjoining the edges of the
support tetrahedron, and four parts of the original hypersurface
adjoining the vertices of the support tetrahedron. Thus in this
case the singular points of the boundary divide it into 15
($1+4+6+4$) strata as shown in Fig.~\ref{f3}, left.

In order to imagine a convex hull in four-dimensional space it is
convenient to project it affinely into the support hyperplane to
the original hypersurface. Then the boundary of the convex hull is
locally the graph of a continuously differentiable function (see,
e.g., \cite{Z}) whose typical singularities are discontinuities of
the second derivative. In a neighborhood of a support tetrahedron
of a generic three-dimensional hypersurface the points of such
discontinuities form 28 smooth strata which are shown in
Fig.~\ref{f3}, left.

According to \cite{S1}, the convex hull of a generic
three-dimensional hypersurface has the simplest singularities at
the interior points of the faces of a support tetrahedron and the
angle singularities at the interior points of its edges. In the
present paper it is proved that the germ of our convex hull at no
vertex of a support tetrahedron contains functional moduli with
respect to diffeomorphisms. This singularity is denoted by ${\cal
R}_3$.  Its normal form is not found nor the precise number of
numerical moduli. It is only proved certainly not to be less than
nine but apparently is much more.

Moreover, we investigate singularities in the only remaining case when a generic
hypersurface has a support $A_1A_3$-plane. The segment between
their tangency points is called the {\it support $A_1A_3$-segment\/}. It lies
entirely in the boundary of the convex hull of the original hypersurface.
According to our results, the convex hull of a generic three-dimensional
hypersurface has the simplest singularities at the interior points of its
support $A_1A_3$-segment and the angle singularity ${\cal R}_2^0$ at the
endpoint $A_3$. At the endpoint $A_1$ of the support $A_1A_3$-segment
there appears one more singularity of the convex hull.  This singularity
is denoted by ${\cal V}_3$, does not contain functional moduli, and is
diffeomorphic to its normal form. This normal form is the under-graph of
the square of the distance to that component of the complement to the
swallowtail which consists of polynomials without real roots.

In a neighborhood of a support $A_1A_3$-segment the boundary of
the convex hull of a generic hypersurface is divided by the
singular points into five strata four of which are smooth and the
fifth one is non-smooth --- this phenomenon has not occurred
before. The projection of the singular points of the boundary of
the convex hull into the support $A_1A_3$-plane is shown in
Fig.~\ref{f3}, right. It consists of three smooth surfaces, the
cut swallowtail, and half of the Whitney umbrella. The support
$A_1A_3$-segment lies in the Whitney umbrella, its endpoint $A_3$
is the starting point of a smooth curve which consists of support
\hbox{$A_3$-points} and is shown by dots in Fig.~\ref{f3}. It
should be noted that in general the normalizing diffeomorphisms
preserve neither the support $A_1A_3$-segment, nor the curve of
support $A_3$-points.

Inside the cut swallowtail the boundary of the convex hull is the
original hypersurface. A smooth stratum webbed from support triangles
which degenerate into the support $A_1A_3$-segment is inside the
Whitney umbrella. The stratum between these two is non-smooth and
webbed from support segments. One more part of the original hypersurface
bounded by two smooth surfaces adjoins the non-singular stratum. The remaining
smooth stratum is bounded by two surfaces too and is webbed from support
segments degenerating into support $A_3$-points.

The union of the cut swallowtail and the half Whitney umbrella
whose intersection lines coincide and whose tangent cones are transversal
is called a {\it sail-boat\/} in \cite{S3}. There it is proved that the
projection of the singular points of the boundary of the convex hull into
a support $A_1A_3$-plane to a generic hypersurface is a sail-boat in a
neighborhood of the point of tangency $A_1$ and all sail-boats are
diffeomorphic to each other. According to our results, in this case the
convex hull itself is reduced to the above local normal form which does not
contain moduli as well.

The following conjecture is formulated in \cite{S2}: the singularities of
the convex hull of a generic compact hypersurface in four-dimensional
affine space do not contain functional moduli. Thus this conjecture is
proved in the present paper.

So in a neighborhood of a typical point of the boundary the convex hull of
a generic compact three-dimensional hypersurface is diffeomorphic to the
closed half-space. The points of the boundary where the convex hull has
the simplest singularities ${\cal R}_1$ form a smooth two-dimensional
surface. Finally, the angle singularities ${\cal R}_2^0$ appear along
smooth curves, the singularities ${\cal R}_2^+$, ${\cal R}_2^-$, ${\cal
R}_3$, and ${\cal V}_3$ appear at isolated points of the boundary of the
convex hull.

The singularities ${\cal R}_1$, ${\cal R}_2^0$, ${\cal R}_2^\pm$, and
${\cal V}_3$ allow the following uniform representations: the simplest
singularity ${\cal R}_1$ is diffeomorphic to the under-graph of the square
of the distance to a three-dimensional half-space, the angle singularity
${\cal R}_2^0$ or ${\cal R}_2^\pm$ is diffeomorphic to the under-graph of
the square of the distance to a dihedral angle of variable value
$\beta(z)$ where $\beta(z)=\beta_0+z$ or $\beta(z)=\beta_0 \pm z^2$
respectively, and the singularity ${\cal V}_3$ is diffeomorphic to the
under-graph of the square of the distance to one of the components of the
complement to the swallowtail. For the singularity ${\cal R}_3$ there is
apparently no analogous representation.

\medskip

{\bf Singularities of Legendre varieties.} It turns out that the
boundary of the convex hull of a generic compact three-dimensional
hypersurface is the front of a Legendre variety which can have
only the following standard singularities: $\widetilde{R_1}$,
$\widetilde{R_2}$, $\widetilde{R_3}$, and $\widetilde{V_3}$. These
singularities do not have moduli with respect to contact
diffeomorphisms. The singularities $\widetilde{R_1}$ form smooth
two-dimensional surfaces, the singularities $\widetilde{R_2}$
appear along smooth curves, and the singularities
$\widetilde{R_3}$ and $\widetilde{V_3}$ appear at isolated points
of the Legendre variety. The singularities $\widetilde{R_1}$ are
projected to the simplest singularities of the convex hull, the
singularities $\widetilde{R_2}$ are projected to the angle ones,
the singularities $\widetilde{R_3}$ are projected to the
singularities ${\cal R}_3$, and the singularities
$\widetilde{V_3}$ are projected to the singularities ${\cal V}_3$.

The singularities $\widetilde{R_l}$ consist of $2^l$ smooth
strata. The singularity $\widetilde{V_3}$ consists of one non-smooth and
two smooth strata. If $l=3$, the family
$$
\Phi_{l}(\tau,\lambda)=\frac{1}{2}(\tau^{l+1}+\lambda_{1}\tau^{l-1}+\dots+
\lambda_{l-1}\tau+\lambda_{l})^2
$$
generates one of the smooth strata of the
singularity $\widetilde{V_3}$ and the non-smooth one. According to
Theorem  8 from \cite{ZR}, all families $\Phi_{l}$
generate Legendre (Lagrangian) varieties
$\widetilde{\Phi_l}$ whose generic
Legendre (Lagrangian) projections are stable in the sense of 3.3
in \cite{G}. Open Whitney umbrellas and open swallowtails possess the
same property as well, but our varieties $\widetilde{\Phi_l}$ are new if
$l\geq 3$ ($\widetilde{\Phi_1}$ is an intersection of curves and
$\widetilde{\Phi_2}$ is the open Whitney umbrella).

\medskip

{\bf Higher dimensions.} According to \cite{S4}, in dimension five
or more the convex hull can have functional moduli which are not
removed by small perturbation of the original hypersurface. For
example, they appear as relations between several numerical moduli
(which are already in four-dimensional space) along the lines
formed by the vertices of the support tetrahedrons.

\medskip

{\bf Terminology.} The term ``smooth'' always means ``infinitely
smooth''. The term ``generic'' is used for smooth mappings and
means that the given proposition is only true for some open
everywhere dense set in the $C^\infty$-space of all mappings
considered.  The term ``typical'' is used for points of varieties
and means that the given proposition is only true for some open
everywhere dense set of points.

\medskip

{\bf Organization of the paper.} In Section \ref{sc1} we give the
necessary definitions and rigorously formulate the results
(Theorems \ref{P} and \ref{A}) outlined in Introduction. In
Section~\ref{sc2} the proof of these results is divided into two
steps which are Theorems \ref{B}, \ref{C}, and \ref{D} proved in
Section \ref{sc3}.

Theorem \ref{B} from Section \ref{sc2} states that the boundary of the
convex hull of a generic compact three-dimensional hypersurface is the
front of a Legendre variety which can have only the standard
singularities: $\widetilde{R_1}$, $\widetilde{R_2}$,
$\widetilde{R_3}$, and $\widetilde{V_3}$. Their normal forms with
respect to contact diffeomorphisms are described in Section
\ref{sc2}.

Theorem \ref{C} from Section \ref{sc2} contains the classification
of germs of a generic Legendre fibration at points of the Legendre
variety up to contact diffeomorphisms preserving its local normal
forms $\widetilde{R_1}$, $\widetilde{R_2}$, $\widetilde{R_3}$, and
$\widetilde{V_3}$ already found. However, we consider only
fibrations with the following property: the front of the Legendre
variety must be a continuously differentiable manifold. (The
boundary of the convex hull of any compact hypersurface is
continuously differentiable.)

\medskip

{\bf Acknowledgments.} The author is grateful to V.~D.~Sedykh for
statement of the problem and interest in the paper, and to
V.~V.~Goryunov and Yu.~G.~Prokhorov for numerous useful
discussions.

\section{Classification of singularities of convex hulls}
\label{sc1}

\subsection*{Convex hulls and support simplexes}

\begin{defin}
The {\it convex hull\/} of a compact subset of an affine space is
the intersection of all closed half-spaces which contain the
subset.
\end{defin}

\medskip

\begin{defin}
A hyperplane is called a {\it support hyperplane\/} to a subset of
an affine space if the subset lies entirely in one of the two
half-spaces defined by the hyperplane and has at least one common
point with the hyperplane.  \end{defin}

\medskip

If the subset is a manifold then the so-called $A_{l_1} \dots
A_{l_m}$-planes are distinguished among its support hyperplanes.

\medskip

\begin{defin}
A support hyperplane to a manifold is called an {\it $A_{l_1}
\dots A_{l_m}$-plane\/} ($l_1$, \dots, $l_m$ are positive odd
numbers) if:

1) it has $m$ points of tangency $A_{l_1}$, \dots, $A_{l_m}$ with the
manifold;

2) these points of tangency are the vertices of an $(m-1)$-dimensional
simplex which is called {\it a support simplex of the  $A_{l_1} \dots
A_{l_m}$-plane\/} or {\it a support $A_{l_1} \dots A_{l_m}$-simplex\/}.

A support $A_1$-plane is called {\it non-singular\/}. The
remaining support hyperplanes are called {\it singular\/}.
\end{defin}

\medskip

\begin{remark}
If $l_1=\dots=l_m=1$ then we write $mA_1$ instead of $A_1 \dots
A_1$.
\end{remark}

\medskip

\begin{remark}
A $k$-dimensional manifold has a point of tangency $A_l$ ($l$ is a
positive integer) with a hypersurface if the restriction of some
local equation of the hypersurface has the form
$\pm\xi_1^{l+1}\pm\xi_2^2\pm\dots\pm\xi_k^2$ in suitable local
coordinates on the manifold.  \end{remark}

\medskip

According to \cite{Z}, if $k=1$ or $n\leq 7$ then the boundary of
the convex hull of a generic compact $k$-dimensional manifold in
$n$-dimensional affine space consists of support $A_{l_1} \dots
A_{l_m}$-simplexes where $l_1+\dots+l_m\leq n$. In particular, in the
case $n=4$ it consists of support $A_1$-points, support
$2A_1$-segments, support $3A_1$-triangles, support $4A_1$-tetrahedrons,
support $A_3$-points, and support $A_1A_3$-segments.

\subsection*{Hierarchy of singularities of convex hulls}

{\bf Simplest singularity ${\cal R}_1$.} A germ of a convex hull
in four-dimensional space with coordinates $(x,y,z,t)$  has the
{\it singularity\/} ${\cal R}_1$ if it is diffeomorphic to the
germ (at the origin) of the set
$$
{\cal R}_1 = \left\{ \min \limits_{p\geq 0} \left(
\frac{p^2}{2}+px+t \right) \leq 0 \right\}.
$$

\medskip

\begin{remark}
The set ${\cal R}_1$ is the under-graph of the function of $x$,
$y$, and $z$ which is equal to half the square of the standard
distance on the real line from the point $x$ to the ray $x \geq
0$.  \end{remark}

\medskip

{\bf Angle singularities ${\cal R}_2^0$ and ${\cal R}_2^\pm$.} A
germ of a convex hull in four-dimensional space with coordinates
$(x,y,z,t)$ has the {\it singularity\/} ${\cal R}_2^0$, ${\cal
R}_2^+$, or ${\cal R}_2^-$ if it is diffeomorphic to the germ (at
the origin) of the set ${\cal R}_2(\alpha)$ where $\alpha(z)=a+z$,
$\alpha(z)=a+z^2$, or $\alpha(z)=a-z^2$ respectively, $|a|<1$, and
$$
{\cal R}_2(\alpha) = {}
$$
$$
{} =\left\{ \min \limits_{p,q\geq 0} \left( \frac{p^2+q^2}{2}+
\alpha(z)pq+px+qy+t \right) \leq 0 \right\}.
$$

\medskip

\begin{remark}
The set ${\cal R}_2(\alpha)$ is the under-graph of the function of
$x$, $y$, and $z$ which is equal to half the square of the
distance from the point $(x,y)$ to the coordinate angle $\{x \geq
0, y \geq 0\}$ in the plane with the Euclidean metric
$$
ds^2= \frac{dx^2-2\alpha(z)dxdy+dy^2}
{1-\alpha^2(z)}.
$$
The value of the angle is $\beta(z)= \pi- \arccos\alpha(z)$.
\end{remark}

\medskip

{\bf Singularity ${\cal R}_3$.} A germ of a convex hull in
four-dimensional space with coordinates $(x,y,z,t)$ has the {\it
singularity\/} ${\cal R}_3$ if it is diffeomorphic to the germ (at
the origin) of the set
$$
{\cal R}_3(F) = \left\{ \min \limits_{p,q,r\geq 0}
F(p,q,r;x,y,z,t) \leq 0 \right\}
$$
where $F$ is a polynomial whose quasihomogeneous expansion has the
form $F=F_2 + F_3+\dots$ \ if $\deg p=\deg q=\deg r=1$, $\deg
x=\deg y=\deg z=1$, and $\deg t=2$,
$$
F_2(p,q,r;x,y,z,t)={}
$$
$$
{}=\frac{p^2+q^2+r^2}{2}+ apq+bpr+cqr+px+qy+rz+t,
$$
and the quadratic form
$(p^2+q^2+r^2)/2+apq+bpr+cqr$ is positive definite.

For example, the set ${\cal R}_3(F_2)$ is the under-graph of half the
square of the distance from the point $(x,y,z)$ to the coordinate angle
$\{x \geq 0, y \geq 0, z \geq 0\}$ in space with the Euclidean metric
defined by the matrix
$$
{\left\| \begin{array}{ccc} 1&a&b\\a&1&c\\b&c&1
\end{array} \right\|}^{-1}.
$$

\medskip

\begin{remark}
The numbers $a$, $b$, and $c$ are moduli of the singularity ${\cal
R}_3$ with respect to diffeomorphisms. These moduli are described
in \cite{S2}. Moreover, Theorem \ref{C} implies that six moduli
are among the coefficients of the quasihomogeneous component $F_3$
of the polynomial $F$.  \end{remark}

\medskip

{\bf Singularity ${\cal V}_3$.} A germ of a convex hull in
four-dimensional space has the {\it singularity\/} ${\cal V}_3$ if
it is diffeomorphic to the germ (at the origin) of the under-graph
${\cal V}_3$ of half the square of the standard distance from the
point $(x,y,z)$ to the body
$$
V_3=\left\{(x,y,z)\in {\R}^3: \forall\tau\in{\R}\ \,
\tau^4+x\tau^2+y\tau+z\geq 0\right\}
$$
which is bounded by the cut
swallowtail and consists of the non-negative polynomials of degree four.

\medskip

{\bf Adjacencies of singularities.} The singularities of convex
hulls adjoin each other in the following way:
$$
\begin{array}{ccccccc} &&&&{\cal R}_2^\pm&&\\ &&&&\downarrow&&\\ {\cal
R}_0&\leftarrow&{\cal R}_1&\leftarrow&{\cal R}_2^0&\leftarrow&{\cal R}_3\\
&&&&\uparrow&&\\ &&&&{\cal V}_3&& \end{array}
$$
where ${\cal R}_0$ denotes the germ of a closed half-space at a
point of its boundary.

\subsection*{Results}

The following theorem is proved in \cite{S1}.

\medskip

\begin{theorem}
\label{P}
The convex hull of a generic compact hypersurface lying in
four-dimen\-sional affine space has the following singularities:

${\cal R}_0$ at the support $A_1$-points and at the interior points of the
support simplexes of the $2A_1$-, $3A_1$-, and $4A_1$-planes;

${\cal R}_1$ at the endpoints of the support $2A_1$-segments, at the
interior points of the sides of the support $3A_1$-triangles, and at the
interior points of the faces of the support $4A_1$-tetrahedrons;

${\cal R}_2^0$ at typical vertices of the support $3A_1$-triangles and at
typical interior points of the edges of the support $4A_1$-tetrahedrons;

${\cal R}_2^\pm$ at the remaining finite number of vertices of the
support $3A_1$-triangles and at the remaining finite number of
interior points of the edges of the support $4A_1$-tetrahedrons.
\end{theorem}

\medskip

Main Theorem \ref{A} of the present paper completes this
classification.

\medskip

\begin{theorem}
\label{A}
The convex hull of a generic compact hypersurface lying in
four-dimen\-sional affine space has the following singularities:

${\cal R}_3$ at the vertices of the support $4A_1$-tetrahedrons;

${\cal R}_1$ at the interior points of the support $A_1A_3$-segments;

${\cal R}_2^0$ at the points of tangency $A_3$ of the support
$A_1A_3$-planes;

${\cal V}_3$ at the points of tangency $A_1$ of the support
$A_1A_3$-planes;

\noindent and the quasidegree of the polynomial $F$ from the
definition of the singularity ${\cal R}_3$ is bounded by a number
$d \geq 3$ which does not depend on the original hypersurface.
\end{theorem}

\medskip

Theorems \ref{P}, \ref{A} follow from Theorems \ref{B}, \ref{C}
formulated in Section \ref{sc2} and proved in Section \ref{sc3} of
the present paper.

\section{Reducing to normal forms}
\label{sc2}

\subsection*{Singularities of Legendre varieties}

Let $B$ be a manifold and $\xi\in B$ be any of its points. A {\it cooriented
contact element\/} to $B$ is any closed linear half-space in the tangent
space $T_\xi B$. The point $\xi$ is called the {\it point of applying\/} of
the contact element.

It is well known that in the space $ST^{*}B$ of all cooriented contact
elements to $B$ there is a natural {\it contact structure\/} (hyperplane
distribution satisfying the condition of maximal non-integrability): a
contact element is allowed to move so that its boundary contains the
velocity of its point of applying. Varieties which are tangent to the
distribution and whose dimension is at most ($\dim B-1$) are called {\it
Legendre\/}.

The hyperplanes of the contact structure in the space $ST^{*}B$ are
naturally {\it cooriented\/} outward: a contact element move in positive
direction if it does not contain the velocity of the point of applying.

Later on we realize $n$-dimensional affine space as an open half-sphere in
$(n+1)$-dimensional Euclidean space and work with the whole sphere which has a
natural projective structure: subspaces passing through the center cut out
planes of various dimensions in the sphere. Each hyperplane divides the
sphere into two half-spheres, and the above definitions of a convex hull
and a support hyperplane are suitable for their subsets.

A cooriented contact element to affine space or sphere is naturally
identified with the pair consisting of the closed half-space and the point
of applying which lies in the boundary of the half-space. Two cooriented
contact elements are called {\it complementary\/} to each other if they
consist of different closed half-spaces having common boundary and the
same point of applying.

\medskip

\begin{defin}
A cooriented contact element is called a {\it support element\/}
to a subset $C$ of an affine space or sphere if it consists of a
closed half-space containing $C$ and a point from $C$. (The point
lies on the boundary of the half-space.) A cooriented contact
element which is complementary to a support element is called an
{\it antisupport element\/}. All support elements to $C$ form a
subset $C^\bot$, all antisupport elements to $C$ form a
subset~$C^\top$.
\end{defin}

\medskip

\begin{defin}
A cooriented contact element is called an {\it infinitesimal
support element\/} to a subset $C$ of a manifold if it is applied
at a point $\xi\in C$ and contains the cone which is tangent to
$C$ at the point $\xi$. A cooriented contact element which is
complementary to an infinitesimal support element is called an
{\it infinitesimal antisupport element\/}. All infinitesimal
antisupport elements to $C$ form a subset $\widetilde{C}$.
\end{defin}

\medskip

\begin{remark}
The sets of all infinitesimal support and antisupport elements (to
a subset of a manifold) are functorial with respect to
diffeomorphisms of the manifold.  In general, this is not true for
the sets of support and antisupport elements. \end{remark}

\medskip

\begin{remark}
Let $\xi\in C\subset B$ be a point of a subset $C$ of a manifold
$B$. Let us consider any Riemannian metric on $B$ and some curve
beginning at the point $\xi$ and possessing the following
property: the distance from a point of the curve to the subset $C$
is an infinitesimal whose degree is more than one if the point
approaches~$\xi$. The cone lying in the tangent space $T_\xi B$
and consisting of the rays which are tangent to all such curves is
called {\it tangent\/} to the subset $C$ at the point $\xi$.
\end{remark}

\medskip

Let $R_0$ denote the hyperplane $s=0$ in the affine space ${\R}^3
\times {\R}$ with coordinates $(u,v,w,s)$. Let us consider the
following subsets of $R_0$:

1) the hyperplane $R_0$ itself;

2) the half-space $R_1=\{(u,v,w)\in R_0: u\geq 0\}$;

3) the dihedral angle $R_2=\{(u,v,w)\in R_0: u\geq 0,\:v\geq 0\}$;

4) the octant $R_3=\{(u,v,w)\in R_0:u\geq 0,\:v\geq 0,\:w\geq 0\}$;

5) the body $V_3=\{(u,v,w)\in R_0: \forall \tau\in{\R}\ \,
\tau^4+u\tau^2+v\tau+w\geq 0\}$ bounded by the cut swallowtail and
consisting of the non-negative polynomials of degree four.

Let $(p,q,r;u,v,w,s)$ be local coordinates on $ST^{*}({\R}^3
\times {\R})$ such that a cooriented contact element applied at
the point $(u,v,w,s)\in {\R}^3 \times {\R}$ has the form
$pdu+qdv+rdw+ds \leq 0$, and let $E_0\subset ST^{*}({\R}^3 \times
{\R})$ denote the space of all such elements. Let
$\widetilde{R_0}$, $\widetilde{R_1}$, $\widetilde{R_2}$,
$\widetilde{R_3}$, and $\widetilde{V_3}\subset E_0$ be the
Legendre varieties consisting of the cooriented contact elements
which are infinitesimal antisupport elements to the subsets $R_0$,
$R_1$, $R_2$, $R_3$, and $V_3 \subset {\R}^3 \times {\R}$
respectively. The variety $\widetilde{R_0}$ is smooth, the
varieties $\widetilde{R_1}$, $\widetilde{R_2}$, $\widetilde{R_3}$,
and $\widetilde{V_3}$ have singularities, for example, at the
origin.

In the coordinates $(p,q,r;u,v,w,s)$:
$$
\begin{array}{l}
\widetilde{R_0}=\{p=q=r=s=0\},\cr \widetilde{R_1}=\{pu=q=r=s=0, \;
p, u\geq 0\},\\
\widetilde{R_2}= \{pu=qv=r=s=0, \; p, u, q, v\geq 0\},\\
\widetilde{R_3}= \{pu=qv=rw=s=0, \; p, u, q, v, r, w\geq 0\}.
\end{array}
$$
The Legendre variety $\widetilde{V_3}$
consists of the following three strata ($\tau$ is a real parameter):
$$
\begin{array}{l}
p=0,\\
q=0,\\
r=0,\\
u\geq -2\tau^2,\\
v=-4\tau^3-2u\tau,\\
w\geq 3\tau^4+u\tau^2,\\
s=0;
\end{array}
\begin{array}{l}
p=r\tau^2,\\
q=r\tau,\\
r\geq 0,\\
u\geq -2\tau^2,\\
v=-4\tau^3-2u\tau,\\
w=3\tau^4+u\tau^2,\\
s=0;
\end{array}
\begin{array}{l}
p=r\tau^2,\\
|q|\leq r|\tau|,\\
r\geq 0,\\
u=-2\tau^2,\\
v=0,\\
w=\tau^4,\\
s=0.
\end{array}
$$
The first and third strata can be extended up to manifolds, the
second stratum can be extended up to an irreducible algebraic
variety. In implicit form this variety is given by the following
polynomials: $32u^3v^2 + 64u^2w^2 + 144uv^2w - 27v^4 + 256w^3$,
\label{pols} $2pu + 3qv + 4rw$, $3pv + 4qw - 2ruv$, $16pw - 8quv -
8ruw - 3rv^2$, $p^2 + qrv + r^2w$, $4pq + r^2v$, $2pr - 2q^2 -
r^2u$, $prv - 4q^2v - 4qrw$, $p^2r - 4pq^2 + r^3w$, and $s$.

Therefore, the Legendre varieties $\widetilde{R_l}$ consist of
$2^l$ strata. They can be extended to manifolds as well as the
first and third strata of $\widetilde{V_3}$. The second stratum of
the Legendre variety $\widetilde{V_3}$ can be extended to an
irreducible algebraic variety. The union of the first and second
strata of $\widetilde{V_3}$ is contact diffeomorphic to the
Legendre variety
$$
\widetilde{\Phi_3}=\{\sigma=\Phi_3(\tau,\lambda),
\Phi_{3,\tau}(\tau,\lambda)=0,\varkappa=-\Phi_{3,\lambda}(\tau,\lambda)\}
$$
(with the contact structure $\varkappa \, d \lambda +d\sigma =0$)
generated by the family
$$
\Phi_{3}(\tau,\lambda)=\frac{1}{2}(\tau^{4}+\lambda_{1}
\tau^{2}+\lambda_{2}\tau+\lambda_{3})^2.
$$
The reducing contact diffeomorphism is given by the formula
$(\varkappa,\lambda,\sigma)=(p,q,r;u,v,\allowbreak w-r,s+r^2/2)$.

\medskip

\begin{defin}
We say that a three-dimensional Legendre variety has the
singularity $\widetilde{R_1}$, $\widetilde{R_2}$,
$\widetilde{R_3}$, or $\widetilde{V_3}$ at some point if its germ
at this point is (up to a local diffeomorphism respecting the
contact structures and their coorientations) respectively the germ
of the Legendre variety $\widetilde{R_1}$, $\widetilde{R_2}$,
$\widetilde{R_3}$, or $\widetilde{V_3}\subset E_0$ at the origin.
\end{defin}

\medskip

These singularities of Legendre varieties adjoin each other in the
following way:
$$
\begin{array}{ccccccc}
\widetilde{R_0}&\leftarrow&\widetilde{R_1}&\leftarrow&\widetilde{R_2}
&\leftarrow&\widetilde{R_3}\\ &          &     &          &\uparrow&
     &     \\ &          &     &          &\widetilde{V_3}   &          &
\end{array}
$$
where the indices are equal to the codimensions of the strata of a Legendre
variety which consist of the corresponding singular points.

\medskip

\begin{theorem}
\label{B}
Let $M\subset S^4$ be a compact smooth hypersurface, $[M]\subset S^4$ be
its convex hull, $[M]^\bot \subset ST^{*}S^4$ be the set of all
cooriented contact support elements to $[M]$, and $\pi:  ST^{*}S^4 \to S^4$ be the
natural projection.

Then $[M]^\bot$ is a Legendre variety which is uniquely projected onto
its front $\pi([M]^\bot)$ which is the boundary of the convex hull of the
original hypersurface $M$. If the hypersurface $M$ is generic, the
Legendre variety $[M]^\bot$ can have only the above singularities
$\widetilde{R_1}$, $\widetilde{R_2}$, $\widetilde{R_3}$, and
$\widetilde{V_3}$. Moreover:

{\rm1)}\enspace the Legendre variety $[M]^\bot$ is smooth above the support
$A_1$-points and above the interior points of the support simplexes of the
$2A_1$-, $3A_1$-, and $4A_1$-planes of the hypersurface $M$;

{\rm2)}\enspace the singularities $\widetilde{R_1}$ appear above the
endpoints of the support $2A_1$-segments, above the interior points of the
sides of the support $3A_1$-triangles, above the interior points of the faces
of the support $4A_1$-tetrahedrons, and above the interior points of the
support $A_1A_3$-segments of the hypersurface $M$;

{\rm3)}\enspace the singularities $\widetilde{R_2}$ appear above the vertices
of the support $3A_1$-triangles, above the interior points of the edges of
the support $4A_1$-tetrahedrons, and above the points of tangency $A_3$ of
the support $A_1A_3$-planes of the hypersurface $M$;

{\rm4)}\enspace the singularities $\widetilde{R_3}$ appear above the vertices
of the support $4A_1$-tetrahedrons of the hypersurface $M$;

{\rm5)} the singularities $\widetilde{V_3}$ appear above the points of
tangency $A_1$ of the support $A_1A_3$-planes of the hypersurface $M$.
\end{theorem}

\subsection*{Normal forms of Legendre fibrations}

A {\it Legendre\/} mapping is a diagram
$$
L^{n-1} \hookrightarrow E^{2n-1} \to B^n
$$
consisting of an embedding of a Legendre variety $L^{n-1}$ into
the space $E^{2n-1}$ with a cooriented contact structure and a
Legendre fibration $E^{2n-1} \to B^n$. (A smooth fibration whose
fibers are Legendre manifolds is called {\it Legendre\/}.) {\it
Equivalence\/} of Legendre mappings is a commutative diagram
$$
\begin{array}{ccccc} L&
\hookrightarrow& E& \to & B\\ \updownarrow&&\updownarrow&&\updownarrow\\
L^\prime& \hookrightarrow& E^\prime& \to & B^\prime \end{array}
$$
where the middle vertical arrow is a diffeomorphism sending the cooriented
contact structures to each other.

Legendre fibrations are locally given with the help of the
generating families defined below.

Let us consider again the space $E_0$ of cooriented contact
elements $\varkappa d\lambda +ds\leq 0$ applied at the points
$(\lambda,s) \in {\R}^3 \times {\R}$ where $\varkappa=(p,q,r)$,
$\lambda=(u,v,w)$. The cooriented contact structure on $E_0$ is
given by the zero subspaces of the form $\varkappa d\lambda +ds$.

Let $F:({\R}^3\times{\R}^3\times{\R},0) \to ({\R},0)$ be a germ of
a family of smooth functions of $\varkappa$ which depend smoothly
on the parameters $\mu=(x,y,z) \in {\R}^3$ and $t \in {\R}$. Let
$F_\varkappa(0)=0$ and let $F$ satisfy the condition of
non-degeneracy at the origin:
$$
\det\left\| \begin{array}{cc} F_{\varkappa\mu}&F_{\varkappa t}\\
F_\mu&F_t \end{array} \right\|\ne 0.
$$
Then $F$ is called the {\it
generating family\/} of the germ
$$
\pi :  (E_0,0) \to ({\R}^3\times{\R},0), \quad \pi
(\varkappa,\lambda,s) =(\mu,t)
$$
of the Legendre fibration whose fibers are given by the formula
$$
\pi^{-1}(\mu,t)= \{(\varkappa,\lambda,s) \in E_0 :
\phantom{AAAAAAAAAAAAA}
$$
$$
\phantom{AAAAAAAA} \lambda=F_\varkappa(\varkappa,\mu,t),
s=F(\varkappa,\mu,t)-\varkappa\lambda\}.
$$
This fibration is correctly defined in a neighborhood of the
origin in view of the non-degeneracy of $F$.

Thus $\varkappa$ are local coordinates on the fibers of the germ $\pi$,
$(\mu,t)$ are local coordinates on its base, and the cooriented contact
structure on $E_0$ is defined by the form $F_\mu d\mu+F_t dt$.

For example, the natural Legendre fibration $(\varkappa,\lambda,s)
\mapsto (\lambda,s)$ is given by the generating family
$\varkappa\mu+t=px+qy+rz+t$.

\medskip

\begin{theorem}
\label{C} Let $L \subset ST^{*}S^4$ be a Legendre variety with the
singularities $\widetilde{R_1}$, $\widetilde{R_2}$,
$\widetilde{R_3}$, or $\widetilde{V_3}$ and $\pi:  ST^{*}S^4 \to
S^4$ be a generic Legendre fibration. If $\pi(L)$ is a
continuously differentiable manifold then the germ of the Legendre
mapping $L \hookrightarrow ST^{*}S^4
\stackrel{\pi}{\longrightarrow} S^4$

{\rm1)}\enspace at each smooth point of the variety $L$ is
equivalent to the germ $(\widetilde{R_0},0) \hookrightarrow
(E_0,0) \to ({\R}^3\times{\R},0)$ of a Legendre mapping such that
its second arrow is given by the generating family
$$
px+qy+rz+t;
$$

{\rm2)}\enspace at each singular point $\widetilde{R_1}$ of the
variety $L$ is equivalent to the germ $(\widetilde{R_1},0)
\hookrightarrow (E_0,0) \to ({\R}^3\times{\R},0)$ of a Legendre
mapping such that its second arrow is given by the generating
family
$$
\frac{p^2}{2}+px+qy+rz+t;
$$

{\rm3)}\enspace at each singular point $\widetilde{R_2}$ of the
variety $L$ is equivalent to the germ $(\widetilde{R_2},0)
\hookrightarrow (E_0,0) \to ({\R}^3\times{\R},0)$ of a Legendre
mapping such that its second arrow is given by one of the
generating families having the form
$$
\frac{p^2+q^2}{2}+\alpha(z)pq+px+qy+rz+t
$$
where $\alpha(z)=a+z$ for a
typical singular point $\widetilde{R_2}$, $\alpha(z)=a\pm z^2$ for the
remaining finite number of singular points $\widetilde{R_2}$, and $|a|<1$ is a
continuous invariant;

{\rm4)}\enspace at each singular point $\widetilde{R_3}$ of the
variety $L$ is equivalent to the germ $(\widetilde{R_3},0)
\hookrightarrow (E_0,0) \to ({\R}^3\times{\R},0)$ of a Legendre
mapping such that its second arrow is given by a polynomial
generating family $F(p,q,r;x,y,z,t)$ whose quasihomogeneous
expansion has the form $F=F_2+\dots+F_d$ if $\deg p=\deg q=\deg
r=1$, $\deg x=\deg y=\deg z=1$, and $\deg t=2$, where $$ F_2=
\frac{p^2+q^2+r^2}{2}+ apq+bpr+cqr+px+qy+rz+t, $$ $a$, $b$, and
$c$ are continuous invariants, the quadratic form
$(p^2+q^2+r^2)/2+apq+bpr+cqr$ is positive definite, among the
coefficients of $F_3$ there are six continuous invariants, and the
number $d \geq 3$ depends on neither $L$ nor $\pi$;

{\rm5)}\enspace at each singular point $\widetilde{V_3}$ of the
variety $L$ is equivalent to the germ $(\widetilde{V_3},0)
\hookrightarrow (E_0,0) \to ({\R}^3\times{\R},0)$ of a Legendre
mapping such that its second arrow is given by the generating
family
$$
\frac{p^2+q^2+r^2}{2}+px+qy+rz+t.
$$
\end{theorem}

\medskip

\begin{remark}
In the case 5) the generating family can be reduced to the simpler
form
$$
\frac{q^2+r^2}{2}+px+qy+rz+t.
$$
\end{remark}

\subsection*{Normal forms of fronts}

The {\it front} of a Legendre mapping $L^{n-1} \hookrightarrow
E^{2n-1} \to B^n$ is the image of $L^{n-1}$ in $B^n$. The
coorientation of the contact structure on $E$ induces a
coorientation of the front that is, in general, a hypersurface
with singularities.

\medskip

\begin{theorem}
\label{D} The cooriented fronts of the germs of Legendre mappings
from the items 1)--5) of Theorem \ref{C} are the germs (at the
origin) of the boundaries of the sets ${\cal R}_0$, ${\cal R}_1$,
${\cal R}_2(\alpha)$, ${\cal R}_3(F)$, and ${\cal V}_3$ which are
cooriented in the direction of the axis $t$.  \end{theorem}

\section{Proofs}
\label{sc3}

\subsection*{Duality}

All closed half-spaces of the sphere $S^n$ form the {\it dual}
sphere $\widehat{S^n}$. Any point of the sphere $S^n$ is
canonically identified with the closure of an open half-space in
the dual sphere $\widehat{S^n}$. This open half-space in
$\widehat{S^n}$ is formed by all closed half-spaces in $S^n$ which
do not contain the given point of the sphere $S^n$. In other
words, the sphere $S^n$ itself is dual to the sphere
$\widehat{S^n}$. Moreover, the spaces of all cooriented contact
elements to the dual spheres $S^n$ and $\widehat{S^n}$ are
canonically identified with each other. The natural cooriented
contact structures coincide after this identification, so later on
we do not distinguish these spaces and denote them uniformly:
$E^{2n-1}=ST^{*} S^n\cong ST^{*}\widehat{S^n}$.

\medskip

\begin{lemma}
\label{L1} Let $M\subset S^n$ be a compact subset, $[M]\subset
S^n$ be its convex hull, and $[M]^\bot \subset E^{2n-1}$ be the
set of all cooriented contact elements which are support to $[M]$.
Let $M^\ast\subset\widehat{S^n}$ denote the set of all closed
half-spaces containing $M$ and $\widetilde{M^\ast} \subset
E^{2n-1}$ denote the set of all cooriented contact elements which
are infinitesimally antisupport to $M^\ast$. Then
$[M]^\bot=\widetilde{M^\ast}$.  \end{lemma}

\medskip

\begin{proof}
Let ${M^\ast}^\top \subset E^{2n-1}$ denote the set of all
cooriented contact elements which are antisupport to $M^\ast$.
Then $[M]^\bot={M^\ast}^\top$. Now let us prove that
${M^\ast}^\top= \widetilde{M^\ast}$.

The inclusion ${M^\ast}^\top\subset \widetilde{M^\ast}$ is
obvious. Let $\mathrm{l}\notin {M^\ast}^\top$ and $\xi$ be the
point of applying of the element $\mathrm{l}$. If $\xi\notin
M^\ast$ then $\mathrm{l}\notin \widetilde{M^\ast}$. If $\xi\in
M^\ast$ then in the interiority of the element $\mathrm{l}$ there
exists a point $\eta\in M^\ast$. But the segment between the
points $\xi$ and $\eta$ lies entirely in $M^\ast$ and,
consequently, in the cone which is tangent to $M^\ast$ at the
point $\xi$. Again $\mathrm{l}\notin \widetilde{M^\ast}$ that
proves the inclusion ${M^\ast}^\top\supset \widetilde{M^\ast}$.
\end{proof}

\subsection*{Proof of Theorem \ref{B}}

Theorem \ref{B} is partially proved in \cite{S1}. It only remains
to investigate the singularities of the Legendre variety
$[M]^\bot$ which appear above the vertices of the support
$4A_1$-tetrahedrons and above the support $A_1A_3$-segments of the
hypersurface $M$.

\medskip

{\bf Case $4A_1$.} Let $O \subset S^4$ be a closed half-space
containing a generic hypersurface $M \subset S^4$ and bounded by a
support to $M$ $4A_1$-plane and $A_1$ be any point of tangency of
the plane and $M$.

According to Lemma \ref{L1}, $[M]^\bot=\widetilde{M^\ast}$. Let us
choose local smooth coordinates $(u,v,w,s)$ on $\widehat{S^4}$
with an origin $O$ such that the set $M^\ast$ is locally given by
the inequalities $u \geq 0$, $v \geq 0$, $w \geq 0$, and $s \geq
0$, and the cooriented contact element $ds \leq 0$ applied at $O$
is the pair $(O,A_1)$. To investigate the singularity of the
Legendre variety $\widetilde{M^\ast}$ above the point $A_1$ let us
consider the cooriented contact elements from $\widetilde{M^\ast}$
which have the form $pdu+qdv+rdw+ds \leq 0$.  The variety of all
such elements is the Legendre variety $\widetilde{R_3}$.

\medskip

{\bf Case $A_1A_3$.} Let $O \subset S^4$ be a closed half-space
containing a generic hypersurface $M \subset S^4$ and bounded by a
support to $M$ $A_1A_3$-plane, $A_1$ and $A_3$ be the points of
tangency ($A_1$ and $A_3$ respectively) of the plane and $M$.

According to Lemma \ref{L1}, $[M]^\bot=\widetilde{M^\ast}$. Let us
choose local smooth coordinates $(u,v,w,s)$ on $\widehat{S^4}$
with an origin $O$ such that the set $M^\ast$ is locally given by
the condition $(u,v,w)\in V_3$ and the inequality $s \geq 0$. Then
the cooriented contact element $ds \leq 0$ applied at $O$ is the
pair $(O,A_1)$, the element $dw \leq 0$ is the pair $(O,A_3)$, and
the elements $dw+r^\prime ds \leq 0$ where $r^\prime \geq 0$
consist of the half-space $O$ and interior points of the support
$A_1A_3$-segment of the hypersurface $M$.

To investigate the singularity of the Legendre variety
$\widetilde{M^\ast}$ above the point $A_1$ let us consider the
cooriented contact elements from $\widetilde{M^\ast}$ which have
the form $pdu+qdv+rdw+ds \leq 0$.  The variety of all such
elements is the Legendre variety $\widetilde{V_3}$.

To investigate the singularities of the Legendre variety
$\widetilde{M^\ast}$ above the point $A_3$ and above the interior
points of the support $A_1A_3$-segment let us consider the
Legendre variety $\widetilde{M^\ast_0}$ of the cooriented contact
elements from $\widetilde{M^\ast}$ which have the form $p^\prime
du+q^\prime dv+dw+r^\prime ds \leq 0$. All such elements are
infinitesimally antisupport to the set $\partial V_3\times {\R}^+$
where $\partial V_3$ is the boundary of the body $V_3$ (the cut
swallowtail) and ${\R}^+$ is the real ray $s \geq 0$. The set
$\partial V_3\times {\R}^+$ has the form $u \geq -2\tau^2$,
$v=-4\tau^3-2u\tau$, $w=3\tau^4+u\tau^2-u^2/4$, $s \geq 0$ where
$\tau$ is a real parameter. The Legendre variety
$\widetilde{M^\ast_0}$ consists of the following four strata:

1) $p^\prime  \geq 0$, $u-2p^\prime +2{q^\prime}^2=0$, $s \geq 0$,
$r^\prime=0$, $v+4p^\prime q^\prime =0$, $w+{p^\prime}^2-4p^\prime
{q^\prime}^2=0$;

2) $p^\prime \geq 0$, $u-2p^\prime+2{q^\prime}^2=0$, $s=0$,
$r^\prime \geq 0$, $v+4p^\prime q^\prime =0$,
$w+{p^\prime}^2-4p^\prime {q^\prime}^2=0$;

3) $p^\prime =0$, $u+2{q^\prime}^2 \leq 0$, $s \geq 0$, $r^\prime
=0$, $v=0$, $w=0$;

4) $p^\prime =0$, $u+2{q^\prime}^2 \leq 0$, $s=0$, $r^\prime \geq
0$, $v=0$, $w=0$.

After the substitution $p^\prime=U$, $q^\prime=W$, $r^\prime=Q$,
$u=-P+2U-2W^2$, $v=-R-4UW$, $w=S+PU+RW-U^2+4UW^2$ the cooriented
contact elements have the form $PdU+QdV+RdW+dS \leq 0$ and the
Legendre variety $\widetilde{M^\ast_0}$ coincides with
$\widetilde{R_2}$ which consists of the following strata:

1) $U\geq 0$, $P=0$, $V\geq 0$, $Q=0$, $R=0$, $S=0$;

2) $U\geq 0$, $P=0$, $V=0$, $Q\geq 0$, $R=0$, $S=0$;

3) $U=0$, $P\geq 0$, $V\geq 0$, $Q=0$, $R=0$, $S=0$;

4) $U=0$, $P\geq 0$, $V=0$, $Q\geq 0$, $R=0$, $S=0$.

\subsection*{Contact vector fields}

A vector field preserving a contact structure on a manifold is
called {\it contact}. If the contact structure is given by the
zero subspaces of a 1-form then any contact vector field is
uniquely defined by its {\it generating function} that is the
substitution of the field into the 1-form.

For example, on the space $E_0$ with the contact structure given
by the form $\varkappa d\lambda+ds$ the contact vector field
defined by a generating function $K(\varkappa,\lambda,s)$ has the
form
\begin{equation} \label{F1} \dot\varkappa=\varkappa K_s-K_\lambda, \quad
\dot\lambda= K_\varkappa, \quad \dot s=K-\varkappa K_\varkappa.
\end{equation} Let us consider the Legendre fibration
$(\varkappa,\lambda,s) \mapsto (\mu,t)$ locally given by a
generating family $F(\varkappa,\mu,t)$. Let $K(F)$ denote the
derivative of the generating family when the Legendre fibration is
acted by the contact vector field with the generating function
$K$.

\medskip

\begin{lemma}
\label{L2} $K(F)=-K(\varkappa,F_\varkappa,F-\varkappa
F_\varkappa)$. \end{lemma}

\medskip

\begin{proof}
Explicit calculation of the action of the vector field (\ref{F1})
on the fibers $\lambda=F_\varkappa$, $s=F-\varkappa F_\varkappa$.
\end{proof}

\medskip

\begin{remark}
The associative algebra of generating functions has a natural
structure of a Lie algebra. Namely, the bracket of two generating
functions is equal to the generating function of the commutator of
the corresponding contact vector fields: $\{K,L\}=\varkappa K_s
L_\varkappa - K_\lambda L_\varkappa  + K_\varkappa L_\lambda + K
L_s - \varkappa K_\varkappa L_s - K_s L$.  \end{remark}

\subsection*{Proof of Theorem \ref{C}}

In the present paper this theorem is proved only for the new
singularities $\widetilde{R_3}$ and $\widetilde{V_3}$. The rest
cases 1)--3) are actually examined in \cite{S1}.

\medskip

{\bf Case $\widetilde{R_3}$.} In a neighborhood of any singular
point $\widetilde{R_3}$ the Legendre variety $L$ consists of eight
strata which can be extended up to manifolds. Because $\pi$ is
generic its restriction on any of the strata is a mapping of
maximal rank.

So in a singular point $\widetilde{R_3}$ of the Legendre variety
$L$ the germ $\pi:(\varkappa,\lambda,s) \mapsto (\mu,t)$, where
$\varkappa=(p,q,r)$, $\lambda=(u,v,w)$, and $\mu=(x,y,z)$, is
given by a generating family $F(\varkappa,\mu,t)$. It follows from
the condition of maximal rank of the restriction of $\pi$ on the
stratum $p=q=r=s=0$. The conditions of maximal rank of the
restriction of $\pi$ on the rest seven strata of $\widetilde{R_3}$
imply that each diagonal minor of the matrix $\left\|
F_{\varkappa\varkappa} \right\|$ is not equal to zero at the
considered singular point $\widetilde{R_3}$.

Moreover, the matrix $\left\|F_{\varkappa\varkappa}\right\|$ is
positive definite at the origin. Indeed, $F_{pp} > 0$, otherwise
the subset $q=r=s=0$ of $\widetilde{R_3}$ consisting of two its
strata is folded. Hence, and $F_{pp}F_{qq}-F_{pq}^2 > 0$,
otherwise the pair of strata $u=r=s=0$ is folded. At last, $\det
\left\| F_{\varkappa\varkappa} \right\| > 0$, otherwise the pair
of strata $u=v=s=0$ is folded.

Let us consider in $F$ the terms of degree one and two with
respect to $\varkappa$, $\mu$, and $t$. Let us reduce them to the
form $F_2$ acting by means of the generating functions $pu$, $qv$,
and $rw$ which preserve $\widetilde{R_3}$, and choosing
coordinates $(\mu,t)$ properly. Now let us fix the degrees of
quasihomogeneity indicated in the formulation of the theorem.

Let ${\cal I}(\widetilde{R_3})$ be the Lie algebra of the germs
(at the origin) of generating functions of contact vector fields
on $E_0$ preserving the origin and the Legendre variety
$\widetilde{R_3}$, and ${\cal D}_{\mu,t}$ be the Lie algebra of
the germs (at the origin) of vector fields on the base of the
Legendre fibration $\pi$ which preserve the origin. Let us
consider the Lie algebra ${\cal G} = {\cal I}(\widetilde{R_3})
\oplus {\cal D}_{\mu,t}$ of the group of the equivalence of
Legendre mappings $(\widetilde{R_3},0) \hookrightarrow (E_0,0) \to
({\R}^3\times{\R},0)$. Generating families of such germs are acted
by the first term of the Lie algebra ${\cal G}$ by means of the
formula from Lemma \ref{L2} and are acted by the second term by
means of derivations. The results of this action of the Lie
algebra ${\cal G}$ are tangent vectors to the space of generating
families. These tangent vectors are elements of the ideal
$$
{\goth p}_{\varkappa,\mu,t} = \left\{ f \in {\cal
E}_{\varkappa,\mu,t} : f \big|_{\varkappa,\mu,t=0} =
f_\varkappa\big|_{\varkappa,\mu,t=0}=0 \right\}
$$
lying in the associative algebra ${\cal E}_{\varkappa,\mu,t}$ of
the germs (at the origin) of smooth functions of
$(\varkappa,\mu,t)$.

\medskip

\begin{lemma}
\label{L3} Let ${\cal S}(F_2)=\{g \in {\cal G} :  g F_2=0\}$ be
the stabilizer of $F_2$. Then there exists a number $N$ with the
following property: if the coefficients $a$, $b$, $c$ and the
cubic quasihomogeneous polynomial $F_3$ are typical then
$\dim_{\R} {\goth p}_{\varkappa,\mu,t}/({\cal G} F_2 + {\cal
S}(F_2)F_3) \leq N$.  \end{lemma}

\medskip

The number $N$ from Lemma \ref{L3} bounds the codimension of the
orbit of the germ $\pi$ in the space of jets of Legendre mappings
with respect to the good geometrical group (\cite{AGLV}, Chapter
3, \S\,2) of the Legendre equivalence preserving the embedding
$\widetilde{R_3} \hookrightarrow E_0$. According to the general
finite-definiteness theorem (\cite{AGLV}, Chapter 3, \S\,2), the
generating family of the germ $\pi$ can be reduced to a polynomial
of quasidegree $d$ depending only on $N$.

The numbers of moduli among the coefficients of the polynomials
$F_2$ and $F_3$ are equal $\dim_{\R} {\goth p}_2 /{\cal G}_0 F_2 =
3$ and $\dim_{\R} {\goth p}_3 /({\cal G}_1 F_2 + {\cal
S}_0(F_2)F_3) = 6$ respectively because ${\cal S}_i (F_2) = 0$ if
$i<0$. Here ${\goth p}_i$, ${\cal G}_i$, and ${\cal S}_i (F_2)$
are quasihomogeneous components of the corresponding algebras of
quasidegree $i$. (See, e.\,g., \cite{AGVI}, \S\,14.)

Thus the considered case $\widetilde{R_3}$ of Theorem \ref{C}
follows from Lemma \ref{L3} which is remained to prove.

\medskip

{\sc Proof of Lemma \ref{L3}.} Let ${\cal E}_{\mu,t}$ be the
associative algebra of the germs (at the origin) of smooth
functions of $(\mu,t)$, ${\goth m}_{\mu,t} \subset {\cal
E}_{\mu,t}$ be the maximal ideal, and
$$
{\goth M} = {\cal
E}_{\varkappa,\mu,t} \langle p \partial_p F_2, q \partial_q F_2, r
\partial_r F_2, F_2 \rangle +{} \phantom{AAAAAA}
$$
$$
\phantom{AAAAAA}{}+ {\goth m}_{\mu,t} \langle p,q,r,1 \rangle +
{\cal E}_{\mu,t} \langle F_3 \rangle \subset {\goth
p}_{\varkappa,\mu,t}
$$
which is a ${\cal E}_{\mu,t}$-submodule where $\partial$ means
partial derivation of the indicated variable.  Let us prove that
${\cal G} F_2 + {\cal S}(F_2)F_3 \supset {\goth M}$.

Indeed, the Lie algebra ${\cal I}(\widetilde{R_3})$ contains the
germs $pu$, $qv$, $rw$, and $s$ which generate an ideal in the
associative algebra of the germs of generating functions. Acting
on $F_2$ this ideal which lies entirely in ${\cal
I}(\widetilde{R_3})$ gives us the first term of ${\goth M}$. The
Lie algebra ${\cal D}_{\mu,t} = {\goth m}_{\mu,t} \langle
\partial_x, \partial_y, \partial_z, \partial_t \rangle$ applied to
$F_2$ gives the second term of ${\goth M}$. At last, ${\cal
S}(F_2)$ is a ${\cal E}_{\mu,t}$-module and $F_3 \in {\cal
S}(F_2)F_3$ because $gF_2=0$ and $gF_3=F_3$ in consequence of
Lemma \ref{L2} where $g=(pu+qv+rw+2s, x\partial_x + y\partial_y +
z\partial_z + 2t\partial_t) \in {\cal I}(\widetilde{R_3}) \oplus
{\cal D}_{\mu,t}$.

Now it remains to prove that $\dim_{\R} {\goth
p}_{\varkappa,\mu,t}/{\goth M} \leq N < \infty$ if $a$, $b$, $c$,
and $F_3$ are typical.

Let us consider ${\cal E}_{\mu,t}$-submodule ${\goth M}^\prime =
{\cal E}_{\varkappa,\mu,t} \langle p \partial_p F_2, q \partial_q
F_2, r \partial_r F_2 \rangle + {\goth m}_{\mu,t} \langle p,q,r,1
\rangle \subset {\goth M}$. If $1-a^2\ne 0$, $1-b^2\ne 0$,
$1-c^2\ne 0$, and $\Delta\ne 0$ where $\Delta =1+2abc-a^2-b^2-c^2$
then ${\cal E}_{\mu,t}$-module ${\goth p}_{\varkappa,\mu,t}/{\goth
M}^\prime$ is generated by the classes $[pqr]$, $[pq]$, $[pr]$,
and $[qr]$ because ${\cal E}_{\mu,t}$-module ${\cal
E}_{\varkappa,\mu,t}/{\cal E}_{\varkappa,\mu,t} \langle p
\partial_p F_2, q \partial_q F_2, r \partial_r F_2 \rangle$
is generated by the monomials $pqr$, $pq$, $pr$, $qr$, $p$, $q$,
$r$, and $1$. The last is implied by the Weierstrass--Malgrange
preparation theorem (see, e.\,g., \cite{AGVI}, Chapter 1, 4.4 and
6.6) applied to the mapping $(\varkappa,\mu,t) \mapsto (\mu,t)$
and the following fact: if the above inequalities are correct then
the monomials generate the linear space ${\cal E}_{\varkappa} /
(p^2+apq+bpr, apq+q^2+cqr, bpr+cqr+r^2)$.

In particular, the following relations are true in ${\goth
p}_{\varkappa,\mu,t}/{\goth M}^\prime$:
$$
\displaylines{[p^2]= [-apq-bpr],\cr [q^2]= [-apq-cqr], \cr [r^2]=
[-bpr-cqr], \cr [p^3]= [-ap^2q-bp^2r-xp^2], \cr [q^3]=
[-apq^2-cq^2r-yq^2], \cr [r^3]= [-bpr^2-cqr^2-zr^2], \cr [p^2q]=
{(1-a^2)}^{-1} [(ac-b)pqr+(ay-x)pq], \cr [pq^2]= {(1-a^2)}^{-1}
[(ab-c)pqr+(ax-y)pq], \cr [p^2r]= {(1-b^2)}^{-1}
[(bc-a)pqr+(bz-x)pr], \cr [pr^2]= {(1-b^2)}^{-1}
[(ab-c)pqr+(bx-z)pr], \cr [q^2r]= {(1-c^2)}^{-1}
[(bc-a)pqr+(cz-y)qr], \cr [qr^2]= {(1-c^2)}^{-1}
[(ac-b)pqr+(cy-z)qr], \cr [p^2qr]= {\Delta}^{-1}
[(-x+c^2x-bcy+ay-acz+bz)pqr], \cr [pq^2r]= {\Delta}^{-1}
[(-bcx+ax-y+b^2y-abz+cz)pqr], \cr [pqr^2]= {\Delta}^{-1}
[(-acx+bx-aby+cy-z+a^2z)pqr].}
$$
Indeed, the expressions for $[p^2]$, $[q^2]$, and $[r^2]$ follow
from $[p\partial_p F_2]=[q\partial_q F_2]=[r\partial_r F_2]=0$,
the expressions for $[p^3]$, $[q^3]$, and $[r^3]$ follow from
$[p^2\partial_p F_2]=[q^2\partial_q F_2]=[r^2\partial_r F_2]=0$,
the expressions for $[p^2q]$ and $[pq^2]$ follow from
$[pq\partial_p F_2]=[pq\partial_q F_2]=0$, the expressions for
$[p^2r]$ and $[pr^2]$ follow from $[pr\partial_p
F_2]=[pr\partial_r F_2]=0$, the expressions for $[q^2r]$ and
$[qr^2]$ follow from $[qr\partial_q F_2]=[qr\partial_r F_2]=0$,
and, at last, the expressions for $[p^2qr]$, $[pq^2r]$, and
$[pqr^2]$ follow from $[pqr\partial_p F_2]=[pqr\partial_q
F_2]=[pqr\partial_r F_2]=0$.

Let us consider in ${\cal E}_{\mu,t} \langle pqr, pq, pr, qr
\rangle$ the ${\cal E}_{\mu,t}$-submodule ${\goth N}$ generated by
the elements $f_2^p$, $f_2^q$, $f_2^r$, $f_2^{qr}$, $f_2^{pr}$,
$f_2^{pq}$, $f_2^{pqr}$, and $f_3 \in {\cal E}_{\mu,t} \langle
pqr, pq, pr, qr \rangle$ defined by the following conditions in
${\goth p}_{\varkappa,\mu,t}/{\goth M}^\prime$:  $[2pF_2]=
[f_2^p]$, $[2qF_2]= [f_2^q]$, $[2rF_2]= [f_2^r]$, $[2pqF_2]=
[f_2^{pq}]$, $[2prF_2]= [f_2^{pr}]$, $[2qrF_2]= [f_2^{qr}]$,
$[2pqrF_2]= [f_2^{pqr}]$, and $[F_3]= [f_3]$. The expressions of
the previous indentation imply that the coordinates of the
generators of the module ${\goth N}$ in the basis $\langle pqr,
pq, pr, qr \rangle$ depend rationally on the numbers $a$, $b$, $c$
and the coefficients of the polynomial $F_3$. For example, if
$a=b=c=0$ we get $f_2^p= ypq+zpr$, $f_2^q= xpq+zqr$, $f_2^r=
xpr+yqr$, $f_2^{pq}= zpqr+(2t-x^2-y^2)pq$, $f_2^{pr}=
ypqr+(2t-x^2-z^2)pr$, $f_2^{qr}= xpqr+(2t-y^2-z^2)qr$, $f_2^{pqr}=
(2t-x^2-y^2-z^2)pqr$.

Let $a=b=c=0$, $F_3=zpq+2ypr+3xqr$, and $N= \dim_{\R} {\cal
E}_{\mu,t} \langle pqr, pq, pr, qr \rangle/{\goth N}$. Then $N<
\infty$. Indeed, let us consider the matrix
($\omega=2t-x^2-y^2-z^2$):
$$
\Omega=\left\| \begin{array}{cccc} 0& 0& z& y\\
0& z& 0& x\\ 0& y& x& 0\\x& \omega + x^2 & 0 & 0\\  y& 0& \omega +
y^2& 0\\ z& 0& 0 & \omega + z^2\\ \omega & 0& 0& 0\\ 0& 3x& 2y& z
\end{array} \right\|
$$
consisting of the coordinates of the generators of the module
${\goth N}$ in the basis $\langle pqr, qr, pr, pq \rangle$. It is
directly verified that the condition ${\rm rank}\:\Omega<4$
implies the relations $x=y=z=t=0$ in ${\C}^4$. So $\dim_{\R} {\cal
E}_{\mu,t}/(M_1, \dots, M_{70}) < \infty$ where $(M_1, \dots,
M_{70})$ is the ideal generated by all $4\times4$-minors of the
matrix $\Omega$. Because all elements $M_ipqr$, $M_ipq$, $M_ipr$,
or $M_iqr$ where $i=1,\dots,70$ lie in ${\goth N}$ we get the
required finite dimensionality.

Now let $a$, $b$, $c$, and $F_3$ be typical. Then $\dim_{\R} {\cal
E}_{\mu,t} \langle pqr, pq, pr, qr \rangle/{\goth N} \leq N$.
Indeed, the set of the coefficients satisfying the inequality is
not empty in consequence of the definition of $N$ and open in the
sense of Zariski (see \cite{H}, 12.7.2).

But $\dim_{\R} {\goth p}_{\varkappa,\mu,t}/{\goth M} \leq
\dim_{\R} {\cal E}_{\mu,t} \langle pqr, pq, pr, qr \rangle/{\goth
N}$ because the homomorphism ${\cal E}_{\mu,t} \langle pqr, pq,
pr, qr \rangle \to {\goth p}_{\varkappa,\mu,t}/{\goth M}$ is
surjective if $a$, $b$, $c$, and $F_3$ are typical and the kernel
of this homomorphism contains ${\goth N}$.  $\Box$

\medskip

{\bf Case $\widetilde{V_3}$.} In a neighborhood of any singular
point $\widetilde{V_3}$ the Legendre variety $L$ consists of three
strata two of which can be extended up to manifolds. The tangent
cone to the third stratum consists of the two Legendre planes
$p=q=r=s=0$ and $p=q=w=s=0$. These equations are obtained from the
quadratic parts $2pu + 3qv + 4rw$, $3pv + 4qw$, $16pw$, $p^2$,
$4pq$, and $2pr - 2q^2$ of the polynomials vanishing on the
stratum (see p.\,\pageref{pols}). Because $\pi$ is generic its
restriction on any of the two smooth strata is a mapping of
maximal rank. Besides, the differential of $\pi$ does not vanish
on the tangent cone to the non-smooth stratum.

So in a singular point $\widetilde{V_3}$ of the Legendre variety
$L$ the germ $\pi:(\varkappa,\lambda,s) \mapsto (\mu,t)$, where
$\varkappa=(p,q,r)$, $\lambda=(u,v,w)$, and $\mu=(x,y,z)$, is
given by a generating family $F(\varkappa,\mu,t)$. It follows from
the condition of maximal rank of the restriction of $\pi$ on the
stratum $p=q=r=s=0$. The condition of maximal rank of the
restriction of $\pi$ on the stratum $p=v=w=s=0$ implies that
$F_{qq}F_{rr}-F_{qr}^2 \ne 0$ at the origin. At last, the Legendre
plane $p=q=w=s=0$ from the tangent cone to the non-smooth stratum
gives the inequality $F_{rr} \ne 0$.

\medskip

\begin{lemma}
\label{L4} Let the germ at the origin of a Legendre fibration
$\pi:(\varkappa,\lambda,s) \mapsto (\mu,t)$ where
$\varkappa=(p,q,r)$, $\lambda=(u,v,w)$, and $\mu=(x,y,z)$ be given
by a generating family $F(\varkappa,\mu,t)$ which satisfies the
following conditions: $F_{qq}F_{rr}-F_{qr}^2 \ne 0$ and $F_{rr}
\ne 0$, and $\pi^\prime:(E_0,0) \to ({\R}^3\times{\R},0)$ be the
germ at the origin of a Legendre fibration which is sufficiently
close to $\pi$. Then the following equivalence of germs of
Legendre mappings takes place:  $$ \begin{array}{ccccc}
(\widetilde{V_3},0)& \hookrightarrow& (E_0,0)&
\stackrel{\pi^\prime} {\longrightarrow} & ({\R}^3\times{\R},0)\\ \updownarrow&&\updownarrow&&\updownarrow\\
(\widetilde{V_3},0)& \hookrightarrow& (E_0,0)& \stackrel{\pi}
{\longrightarrow} & ({\R}^3\times{\R},0) \end{array}$$ where the
vertical arrows are close to the identity mappings.
\end{lemma}

\medskip

The condition of non-degeneracy divides the space of the
generating families into two connected components transferred to
each other by the substitution $t \mapsto -t$. Each of them is
divided into four components by the conditions
$F_{qq}F_{rr}-F_{qr}^2 \ne 0$ and $F_{rr} \ne 0$. According to
Lemma \ref{L4}, the germs of Legendre mappings are equivalent if
they are given by generating families from the same component.
Hence, each of them can be given by the generating family $(p^2
\pm q^2 \pm r^2)/2 +px+qy+rz+t$ or even by $(\pm q^2 \pm r^2)/2
+px+qy+rz+t$ (see Note to Theorem \ref{C}). It remains to check
that the image $\pi(L)$ is not a continuously differentiable
manifold if there is at least one minus in these families.

Thus the considered case $\widetilde{V_3}$ of Theorem \ref{C}
follows from Lemma \ref{L4} which is remained to prove.

\medskip

{\sc Proof of Lemma \ref{L4}.} The variety $\widetilde{V_3}$ can
be replaced by its extension up to the irreducible algebraic
Legendre variety $\overline{V_3}$ which consists of the following
three strata:  $\{p=q=r=s=0\}$, $\{p=r\tau^2+ru/2, q=r\tau,
v=-4\tau^3-2u\tau, w=3\tau^4+u\tau^2-u^2/4, s=0\}$ ($\tau$ is a
parameter), and $\{p=v=w=s=0\}$. Indeed, any diffeomorphism being
sufficiently close to identity one and preserving $\overline{V_3}$
preserves $\widetilde{V_3}$ as well.

Let us consider the algebra ${\cal E}_\varkappa$ of the germs at
the origin of smooth functions on the fiber $\pi^{-1}(0)$ and the
algebra ${\cal Q}$ of the germs at the origin of smooth functions
on the intersection $\overline{V_3}\cap\pi^{-1}(0)$. There is a
natural restriction homomorphism $i : {\cal E}_\varkappa \to {\cal
Q}$. Let us show that ${\goth m}_\varkappa^2 \subset \ker i$ where
${\goth m}_\varkappa \subset {\cal E}_\varkappa$ is the maximal
ideal. Indeed, $\pi^{-1}(0)= \{(\varkappa,\lambda,s) \in E_0 :
\lambda=f_\varkappa(\varkappa), s=f(\varkappa)- \varkappa
f_\varkappa(\varkappa) \}$ where $f(\varkappa)=F(\varkappa,0,0)$
and the functions $2pu+3qv+4rw$, $3pv+4qw-2ruv$,
$16pw-8quv-8ruw-3rv^2$, $p^2+qrv+r^2w$, $4pq+r^2v$, and $-2s$
vanish on $\overline{V_3}$. So the germs $2pf_p+3qf_q+4rf_r$,
$3pf_q+4qf_r-2rf_pf_q$, $16pf_r-8qf_pf_q-8rf_pf_r-3rf_q^2$,
$p^2+qrf_q+r^2f_r$, $4pq+r^2f_q$, and $2pf_p+2qf_q+2rf_r-2f$ lie
in the kernel of $i$. All of them are in ${\goth m}_\varkappa^2$
and their expansions up to ${\goth m}_\varkappa^3$ in the basis
$\langle p^2, pq, pr, q^2, qr, r^2 \rangle$ form the matrix
$$\left\| \begin{array}{cccccc} \ast& \ast& \ast& 3f_{qq}(0)& 7f_{qr}(0)&
4f_{rr}(0)\\ \ast& \ast& \ast& 4f_{qr}(0)& 4f_{rr}(0)& 0\\ \ast&
\ast& 16f_{rr}(0)& 0& 0& 0\\ 1& 0& 0& 0& 0& 0\\ 0& 4& 0& 0& 0& 0\\
\ast& \ast& \ast& f_{qq}(0)& 2f_{qr}(0)& f_{rr}(0) \end{array}
\right\|$$ which is non-degenerate in accordance with the
condition of the lemma. Hence, the germs generate ${\goth
m}_\varkappa^2/{\goth m}_\varkappa^3$ as linear space over ${\R}$.
According to the Nakayama lemma, they generate ${\goth
m}_\varkappa^2$ as a ${\cal E}_\varkappa$-module. Therefore
${\goth m}_\varkappa^2 \subset \ker i$.

So the germ $(\overline{V_3},0) \hookrightarrow (E_0,0)
\stackrel{\pi}{\longrightarrow} ({\R}^3\times{\R},0)$ of a
Legendre mapping is versal in the sense of 3.3 in \cite{G} because
${\goth m}_\varkappa^2 \subset \ker i$. According to Theorem
$3^\prime$ from \cite{G}, the germ is stable with respect to
perturbations of $\pi$. It remains to note that the point of
applying of the germ is preserved by its equivalence.  $\Box$

\subsection*{Proof of Theorem \ref{D}}

According to the definition, if a germ of a Legendre fibration is
given by a generating family $F(p,q,r;u,v,w,s)$ then its fibers
are defined by the equations $u=F_p$, $v=F_q$, $w=F_r$,
$s=F-pF_p-qF_q-rF_r$ and the cooriented contact structure is
defined by the form $F_x dx+F_y dy+F_z dz+F_t dt$.

The case 1) of the variety $\widetilde{R_0}$ is trivially checked
explicitly.

In the cases 2)--4) substituting the above equations of the fibers
into the explicit coordinate expressions for the varieties
$\widetilde{R_1}$, $\widetilde{R_2}$, and $\widetilde{R_3}$ we get
that the minimum of $F$ is equal to zero as is required.

In the case 5) the equations of the fibers have the form $u=p+x$,
$v=q+y$, $w=r+z$, $s=-(p^2+q^2+r^2)/2+t$ and for the Legendre
fibration itself we get the formula $(p,q,r;u,v,w,s) \mapsto
(u-p,v-q,w-r,s+(p^2+q^2+r^2)/2)$. If $(p,q,r;u,v,w,s) \in
\widetilde{V_3}$ then $s=0$. If in addition $(p,q,r)=0$ then
$(u,v,w) \in V_3$. If $(p,q,r) \ne 0$ then the cooriented contact
element $pdu+qdv+rdw \leq 0$ applied at the point $(u,v,w)$ is
infinitesimally antisupport to $V_3 \subset {\R}^3$. In other
words, $(u,v,w) \in \partial V_3$ and the vector $(-p,-q,-r)$ is
perpendicular to the boundary $\partial V_3$ of the body $V_3$ and
is directed outside.

\bibliography{hull,shock}

\def\cprime{$'$} \def\cprime{$'$}
\begin{thebibliography}{10}

\bibitem{AGLV}
V.~I. Arnol{\cprime}d, V.~V. Goryunov, O.~V. Lyashko, and V.~A. Vassiliev.
\newblock {\em Dynamical systems VI}, volume~6 of {\em Encyclopaedia of
  Mathematical Sciences}.
\newblock Springer-Verlag, 1993.

\bibitem{AGVI}
V.~I. Arnol{\cprime}d, S.~M. Guse{\u\i}n-Zade, and A.~N. Varchenko.
\newblock {\em Singularities of differentiable maps. {I}}.
\newblock Birkh\"auser Boston Inc., Boston, MA, 1985.

\bibitem{G}
A.~B. Givental$'$.
\newblock Singular {L}agrangian manifolds and their {L}agrangian mappings.
\newblock {\em J.~Soviet Math.}, 52:3246--3278, 1990.

\bibitem{H}
R.~Hartshorne.
\newblock {\em Algebraic geometry}, volume~52 of {\em Graduate Texts in
  Mathematics}.
\newblock Springer-Verlag, 1977.

\bibitem{S2}
V.~D. Sedykh.
\newblock Singularities of convex hulls.
\newblock {\em Siberian Math. J.}, 24:447--461, 1983.

\bibitem{S4}
V.~D. Sedykh.
\newblock Functional moduli of singularities of convex hulls of manifolds of
  codimension $1$ and $2$.
\newblock {\em Math. USSR-Sb.}, 47:223--236, 1984.

\bibitem{S1}
V.~D. Sedykh.
\newblock Stabilization of the singularities of convex hulls.
\newblock {\em Math. USSR-Sb.}, 63:499--505, 1989.

\bibitem{S3}
V.~D. Sedykh.
\newblock The sewing of the swallowtail and the {W}hitney umbrella in a
  four-dimensional control system.
\newblock {\em Trudy GANG im. I.~M.~Gubkina}, pages 58--68, 1997.

\bibitem{Z}
V.~M. Zakalyukin.
\newblock Singularities of convex hulls of smooth manifolds.
\newblock {\em Functional Anal. Appl.}, 11:225--227, 1978.

\bibitem{ZR}
V.~M. Zakalyukin and R.~M. Roberts.
\newblock Stability of {L}agrangian manifolds with singularities.
\newblock {\em Functional Anal. Appl.}, 26:174--178, 1992.

\end{thebibliography}
\bibliographystyle{plain}

\end{document}